\documentclass[10pt]{amsart}

\usepackage[all]{xypic}

\usepackage{latexsym}

\usepackage{amssymb}
\usepackage{amsfonts}
\usepackage{amscd}
\usepackage{amsmath,amsthm}

\newtheorem{lemma}{Lemma}[section]

\newtheorem{lem}[lemma]{Lemma}
\newtheorem{prop}[lemma]{Proposition}
\newtheorem{thm}[lemma]{Theorem}
\newtheorem{cor}[lemma]{Corollary}

{

}

\theoremstyle{definition}

\theoremstyle{remark}

\numberwithin{equation}{section}

\newenvironment{pf}{\noindent{\bf Proof.}}{\hfill $\square$\medskip}

\def\PP{{\mathbb P}}

\def\ZZ{{\mathbb Z}}

\def\0ol{{\bar 0}}
\def\1ol{{\bar 1}}
\def\2ol{{\bar 2}}
\def\ol2{{\bar 2}}
\def\3ol{{\bar 3}}
\def\4ol{{\bar 4}}
\def\5ol{{\bar 5}}
\def\6ol{{\bar 6}}
\def\7ol{{\bar 7}}
\def\8ol{{\bar 8}}
\def\9ol{{\bar 9}}

\def\bold0{{\bf 0}}
\def\bold1{{\bf 1}}
\def\bold2{{\bf 2}} 
\def\bold3{{\bf  3}}
\def\bold4{{\bf 4}}
\def\bold5{{\bf 5}}
\def\bold6{{\bf 6}}
\def\bold7{{\bf 7}}
\def\bold8{{\bf 8}}
\def\bold9{{\bf 9}}

\def\P2Skly{\PP^2_{Skly}}

\def\gr{\operatorname {gr}}

\def\th{\operatorname {th}}    

\def\Ann{\operatorname{Ann}}
\def\Aut{\operatorname{Aut}}

\def\fchar{\operatorname{char}}

\def\det{\operatorname{det}}

\def\dim{\operatorname{dim}}

\def\Fdim{{\sf Fdim}}
\def\fdim{{\sf fdim}}

\def\Gr{{\sf Gr}}

\def\liminj{\varinjlim}

\def\mod{{\sf mod}}
\def\Mod{{\sf Mod}}

\def\Proj{\operatorname{Proj}}

\def\QGr{\operatorname{\sf QGr}}
\def\qgr{\operatorname{\sf qgr}}

\def\Span{\operatorname{Span}}

\def\ul1{\operatorname{\underline{1}}}

\def\l{\leftarrow}

\def\a{\alpha}
\def\b{\beta}
\def\c{\gamma}

\def\l{\lambda}

\def\ve{\varepsilon}

\def\cal{\mathcal}

\def\cO{{\cal O}}

\def\Qcoh{{\sf Qcoh}}

\def\dirlim{\mathop{\vtop{\baselineskip -100pt\lineskip -1pt\lineskiplimit 0pt
\setbox0\hbox{lim}\copy0\hbox to \wd0{\rightarrowfill}}}\limits}
\def\invlim{\mathop{\vtop{\baselineskip -100pt\lineskip -1pt\lineskiplimit 0pt
\setbox0\hbox{lim}\copy0\hbox to \wd0{\leftarrowfill}}}\limits}

\def\I11{{1 \kern -0.8pt \! \mbox{l}}}
\def\mumu{{\mu\kern-4.2pt\mu}}
\def\bfmu{{\mu\kern-4.2pt\mu}}
\def\2slash{\backslash \! \backslash}

\def\boxtimes{\setbox0\hbox{$\Box$}\copy0\kern-\wd0\hbox{$\times$}}

\def\hdot{{\:\raisebox{2pt}{\text{\circle*{1.5}}}}}

\begin{document}

\title[Degenerate 3-dimensional Sklyanin algebras]{``Degenerate'' 3-dimensional Sklyanin algebras are monomial algebras}

\author{S. Paul Smith}
\address{ Department of Mathematics, Box 354350, Univ.  Washington, Seattle, WA 98195}
\email{smith@@math.washington.edu}

\maketitle

\begin{abstract}
The 3-dimensional Sklyanin algebras, $S_{a,b,c}$, form a flat family parametrized by points 
$(a,b,c) \in \PP^2-{\frak D}$
where ${\frak D}$ is a set of 12 points. When $(a,b,c)\in {\frak D}$, the algebras having the same defining relations as the  3-dimensional Sklyanin algebras are called ``degenerate Sklyanin algebras''. C. Walton showed they do not have the same properties as the non-degenerate ones. Here we prove that a degenerate Sklyanin algebra is isomorphic to the free algebra on $u$, $v$, and $w$, modulo either the relations $u^2=v^2=w^2=0$ or the relations $uv=vw=wu=0$.  
These monomial algebras are Zhang twists of each other. 
Therefore all degenerate Sklyanin algebras have the same category of graded modules. 
A number of properties of the degenerate Sklyanin algebras follow from this observation.  
We exhibit a quiver $Q$ and an ultramatricial algebra $R$ such that if $S$ is a degenerate Sklyanin algebra then the categories $\QGr S$, $\QGr kQ$, and $\Mod R$, are equivalent.  Here $\QGr(-)$ denotes the category of graded right modules modulo the full subcategory of graded modules that are the sum of their finite dimensional submodules. The group  of cube roots of unity, $\mu_3$, acts as automorphisms 
of the free algebra on two variables, $F$, in such a way that $\QGr S$ is equivalent to $\QGr(F \rtimes \mu_3)$. 
\end{abstract}

\section{Introduction}

Let $k$ be a field having a primitive cube root of unity $\omega$.

\subsection{}
Let ${\frak D}$ be the subset of the projective plane $\PP^2_k$ consisting of the 12 points:
$$
{\frak D}:=\big\{ (1,0,0), \, (0,1,0),\, (0,0,1) \big\} \, \sqcup \, \big\{ (a,b,c) \; | \; a^3=b^3=c^3 \big\}.
$$

The points $(a,b,c) \in \PP^2 - {\frak D}$ parametrize the 3-dimensional Sklyanin algebras,  
$$
S_{a,b,c}=\frac{k\langle x,y,z\rangle}{(f_1,f_2,f_3)}\, ,
$$
where
\begin{align*}
f_1=&ayz+bzy+cx^2  \\
f_2=&azx+bxz+cy^2  \\
f_3=& axy+byx+cz^2.
\end{align*}
In the late 1980s, Artin, Tate, and Van den Bergh,  \cite{ATV1} and \cite{ATV2}, showed that $S_{a,b,c}$ behaves much like the commutative polynomial ring on 3 indeterminates. In many ways it is more interesting because its fine structure is governed by an elliptic curve endowed with a translation automorphism.  
 
 \subsection{}
When $(a,b,c) \in {\frak D}$ we continue to write $S_{a,b,c}$ for the algebra with the same generators and relations and call it a {\sf degenerate} 3-dimensional Sklyanin algebra. This paper concerns their structure.
 
Walton \cite{W} has shown that the degenerate Sklyanin algebras are nothing like the others:
if  $(a,b,c) \in {\frak D}$, then $S_{a,b,c}$ has infinite global dimension, is not noetherian, has exponential growth, and has zero divisors,  none of which happens when $(a,b,c) \notin {\frak D}$.

Unexplained terminology in this paper can be found in Walton's paper \cite{W}.  
 
 \subsection{Results}
  
  Our results show that the degenerate Sklyanin algebras are rather well-behaved, albeit on their own terms.
 
 \begin{thm}
\label{thm.qgr1}
Let $S=S_{a,b,c}$ be a degenerate Sklyanin algebra. 
\begin{enumerate}
  \item 
 If $a=b$, then $S$ is isomorphic to $k\langle u,v,w \rangle$ modulo $u^2=v^2=w^2=0$.
  \item 
If $a\ne b$, then $S$ is isomorphic to $k\langle u,v,w \rangle$ modulo  $uv=vw=wu=0$.
\end{enumerate}
\end{thm}

\begin{cor}
If $(a,b,c)$ and $(a',b',c')$ belong to ${\frak D}$, then $S_{a,b,c}$ and $S_{a',b',c'}$ are Zhang twists of one another, and their categories of graded right modules are equivalent,
$$
\Gr S_{a,b,c} \equiv \Gr S_{a',b',c'}.
$$
\end{cor}

Of more interest to us is a quotient category, $\QGr S$, of the category of graded modules. 
This has surprisingly good properties. We now state the result and then define  $\QGr S$.

\begin{thm}
\label{thm.qgr2}
Let $S=S_{a,b,c}$ be a degenerate Sklyanin algebra. 
There is a quiver $Q$ and an ultramatricial algebra $R$, both independent of $(a,b,c) \in {\frak D}$,
 such that
$$
\QGr S_{a,b,c} \equiv \QGr kQ \equiv \Mod R 
$$
where the algebra $R$ is a direct limit of algebras $R_n$, each of which is a product of three matrix algebras 
  over $k$, and hence a von Neumann regular ring.
  
   Furthermore, there is an action of $\mu_3$, the cube roots of unity in $k$, as automorphisms of
    the free algebra 
  $F=k \langle X,Y \rangle$ such that 
  $$
  \QGr S_{a,b,c} \equiv \QGr(F \rtimes \mu_3).
  $$
\end{thm}

\subsection{}

We work with {\it right modules} throughout this paper.

Let $A$ be a connected graded $k$-algebra.
If $\Gr A$ denotes the category of $\ZZ$-graded right $A$-modules and $\Fdim A$ is the full subcategory of modules that are direct limits of their finite dimensional submodules, then the quotient category
$$
\QGr A:=\frac{\Gr A}{\Fdim A}
$$
plays the role of the quasi-coherent sheaves on a ``non-commutative scheme'' that  we call
$\Proj_{nc} A$.

Artin, Tate, and Van den Bergh, showed that if $(a,b,c) \in \PP^2-{\frak D}$, then $\QGr S$ has ``all''
the properties enjoyed by $\Qcoh \PP^2$, the category of quasi-coherent sheaves on the projective plane
\cite{ATV1}, \cite{ATV2}.

\subsection{Consequences}
 The ring $R$ in Theorem \ref{thm.qgr2} is von Neumann regular because each $R_n$ is. 
 The global dimension of $R$ is therefore equal to 1.

Let $S=S_{a,b,c}$ be a degenerate Sklyanin algebra. Finitely presented monomial algebras 
are coherent so $S$ is coherent. The full subcategory $\gr S$ of $\Gr S$ consisting of the 
finitely presented graded modules is therefore abelian. We write $\fdim S$ for the  full subcategory $\Gr S$ 
consisting of the finite dimensional modules; we have $\fdim S = ( \gr S) \cap (\Fdim S)$. Since $\fdim S$ is a Serre subcategory of $\gr S$ we may form the quotient category 
$$
\qgr S:=\frac{\gr S}{\fdim S}.
$$
The equivalence $ \QGr S \equiv \Mod R$ restricts to an equivalence
$$
\qgr S \equiv \mod R
$$
where $\mod R$ consists of the finitely presented $R$-modules. 

Modules over von Neumann regular rings are flat so finitely presented $R$-modules are projective,
whence the next result. 

\begin{cor}
Every object in $\qgr S$ is projective in $\QGr S$ and every short exact sequence in $\qgr S$ splits. 
\end{cor}

The corollary illustrates, again, that the degenerate Sklyanin algebras are unlike the non-degenerate ones
but still ``nice''.  

Since $R$ is only determined up to Morita equivalence 
there are different Bratteli diagrams 
corresponding to
different algebras in the Morita equivalence class of $R$. There are, however, two ``simplest'' ones,
namely the stationary Bratteli diagrams  that begin
$$
 \UseComputerModernTips
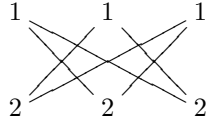
\xymatrix{ 
1 \ar@{-}[drr] \ar@{-}[dr] &1\ar@{-}[dl] \ar@{-}[dr] & 1\ar@{-}[dl] \ar@{-}[dll]
\\
2  & 2  & 2 
}
$$
and
$$
 \UseComputerModernTips
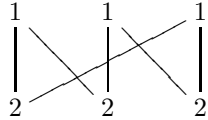
\xymatrix{ 
1 \ar@{-}[d] \ar@{-}[dr] & 1 \ar@{-}[d] \ar@{-}[dr] & 1 \ar@{-}[d] \ar@{-}[dll]
\\
2  & 2  & 2 
}
$$
These diagrams correspond to the quivers $Q$ and $Q'$ that  appear in section  \ref{sect.2.Q}.

\subsection{The Grothendieck group $K_0(\qgr S)$}
Let $S$ be a degenerate Sklyanin algebra and $K_0(\qgr S)$ the Grothendieck group of finitely generated projectives in $\qgr S$.

We write $\cO$ for $S$ as an object in $\qgr S$. The left ideal of $A$ generated by $u-v$ and $v-w$ is free and
$A/A(u-v) + A/(v-w)$ is spanned by the images of 1 and $u$, so $\cO \cong \cO(-1) \oplus \cO(-1)$. 
We show that 
$$
K_0(\qgr S)=  \ZZ\bigg[\frac{1}{8}\bigg] \oplus \ZZ \oplus \ZZ
$$
with $[\cO]=(1,0,0)$.

\subsection{The word ``degenerate''}

The algebras $S_{a,b,c}$ form a flat family on $\PP^2-{\frak D}$. The family does not extend to a flat family on
$\PP^2$ because the Hilbert series for a degenerate Sklyanin algebra is different from the Hilbert series of
the non-degenerate ones. This is the reason for the quotation marks around ``degenerate''. It is not unreasonable to think that there might be another compactification of $\PP^2-{\frak D}$ that parametrizes 
a flat family of algebras that are the Sklyanin algebras on  $\PP^2-{\frak D}$. The obvious candidate to consider is $\PP^2$ blown up at ${\frak D}$.

\subsection{Acknowledgements}

I thank Chelsea Walton for carefully reading this paper and for several useful conversations about its contents.

\section{The degenerate Sklyanin algebras are monomial algebras}

Let
\begin{equation}
\label{eq.AA}
A=\frac{k \langle u,v,w\rangle}{(u^2,v^2,w^2)}
\qquad \hbox{and} \qquad 
A'=\frac{k \langle u,v,w\rangle}{(uv,vw,wu)}.
\end{equation}

The next result uses the blanket hypothesis that the base field $k$ contains a primitive cube root of 
unity which implies that $\fchar k \ne 3$.

\begin{thm}
\label{thm.monomial}
Suppose $(a,b,c) \in {\frak D}$. 
Then
$$
S_{a,b,c} \cong 
\begin{cases}
	A & \text{if $a=b$}
	\\
	A' & \text{if $a \ne b$.}
\end{cases}
$$	
\end{thm}
\begin{pf}
The result is a triviality when $(a,b,c) \in \big\{ (1,0,0), \, (0,1,0),\, (0,0,1) \big\}$
so we assume that $a^3=b^3=c^3$ for the rest of the proof. 

If $\l$ is a non-zero scalar, then $S_{a,b,c} \cong S_{\l a,\l b, \l c}$ so it suffices to prove the theorem 
when 
\begin{enumerate}
  \item 
  $(a,b,c)=(1,1,1)$,
  \item 
  $(a,b,c)=(1,1,c)$ with $c^3=1$ but $c \ne 1$, and
  \item 
  $a \ne b$ and $abc\ne 0$. 
\end{enumerate}
We consider the three cases separately. 

(1) 
Suppose   $(a,b,c)=(1,1,1)$. 
Let $\omega$ be a primitive cube root of unity. Then
\begin{align*}
(x+y+z)^2= & f_1+f_2+f_3,
\\
(x+\omega y+\omega^2 z)^2= & f_1+\omega^2 f_2+\omega f_3,
\\
 (x+\omega^2 y+\omega z)^2 = &  f_1+\omega f_2+\omega^2 f_3.
\end{align*}
Since $\{(1,1,1),(1,\omega,\omega^2), (1,\omega^2,\omega)\}$ is linearly independent,
$$
{\rm span}\{x+y+z, \, x+\omega y+\omega^2 z, \, x+\omega^2 y+\omega z\}=kx+ky+kz.
$$
Therefore $S_{1,1,1} \cong A$.

(2)
Suppose $c \ne 1$ and $(a,b,c)=(1,1,c)$. Then 
\begin{align*}
(x+y+c^{-1}z)^2= & c^{-1} f_1+c^{-1} f_2+f_3,
\\
(x+c^{-1} y+ z)^2= &c^{-1} f_1+ f_2+c^{-1} f_3,
\\
 (c^{-1}x+y+ z)^2 = &  f_1+ c^{-1} f_2+c^{-1} f_3.
\end{align*}
Since 
$$
\det
\begin{pmatrix}
  1 & 1    & c^{-1}   \\
  1 & c^{-1} & 1 \\
  c^{-1} & 1 & 1  
\end{pmatrix} = (c^{-1}-1)^2(c^{-1}+2) \ne 0
$$
$\{x+y+c^{-1}z, \, x+c^{-1} y+  z, \, c^{-1} x+  y+ z\}$ is linearly independent. Hence $S_{1,1,c} \cong A$.

(3)
Suppose $a \ne b$. Let
\begin{align*}
u & = a^{-1}x+b^{-1}y + c^{-1} z,
\\
v & = b^{-1}x+a^{-1}y + c^{-1} z,
\\
w & = abc(x+y)+z.
\end{align*}
Because $\fchar k \ne 3$, the hypothesis that $a \ne b$, implies $\{u,v,w\}$ is linearly independent.
Furthermore,
\begin{align*}
uv & = (abc)^{-1}(f_1+f_2) +f_3  
\\
vw & = af_1+bf_2+cf_3,
\\
wu & = bf_1+af_2+cf_3.
\end{align*}
Hence $S_{a,b,c} \cong A'$.
\end{pf}

\begin{prop}
\label{prop.Z.twist}
The algebras $A$ and $A'$ are Zhang twists of each other.
\end{prop}
\begin{pf} 
The defining relations of $A'$ are $uv=vw=wu=0$ so the map 
$$
\tau(u)=v, \quad \tau(v)=w, \quad \tau(w) =u
$$
extends to an algebra automorphism of $A'$. The Zhang twist of $A'$ by $\tau$ is therefore the 
algebra generated  by $u,v,w$ with defining relations
$$
u*u=u\tau(u)=uv=0, \quad v*v=v\tau(v)=vw=0, \quad w*w=w\tau(w)=wu=0.
$$
The Zhang twist of $A'$ by $\tau$ is therefore isomorphic to $A$.

The defining relations of $A$ are $u^2=v^2=w^2=0$ so the map 
$$
\theta(u)=w, \quad \theta(v)=u, \quad \theta(w) =v
$$
extends to an algebra automorphism of $A$. The Zhang twist of $A$ by $\theta$ is therefore the 
algebra generated  by $u,v,w$ with defining relations
$$
u*v=u\theta(v)=u^2=0, \quad v*w=v\theta(w)=v^2=0, \quad w*u=w\theta(u)=w^2=0.
$$
The Zhang twist of $A$ by $\theta$ is therefore isomorphic to $A'$.
\end{pf}

  Zhang \cite{Z} proved that if a connected graded algebra generated in degree one is a Zhang twist of another one, then their graded module categories are equivalent. This, and Proposition \ref{prop.Z.twist},
 implies the next result.

\begin{cor}
\label{cor.qgr=}
Let $S$ be a degenerate Sklyanin algebra. There are category equivalences
$$
\Gr S \equiv \Gr A \equiv \Gr A'
$$
and
$$
\QGr S \equiv \QGr A \equiv \QGr A'.
$$
In particular, $\Gr S_{a,b,c}$, and hence $\QGr S_{a,b,c}$, is the same for all $(a,b,c) \in {\frak D}$. 
\end{cor}

Walton \cite{W} determined the Hilbert series of $S$ by exhibiting $S$ as a free module of rank 2 over a free subalgebra. Here is an alternative derivation of the Hilbert series using the description of $S$ as a monomial
algebra.

\begin{cor}
\label{cor.H.series}
Let $(a,b,c) \in \frak{D}$. The Hilbert series of $S_{a,b,c}$ is $(1+t)(1-2t)^{-1}$.
\end{cor}
\begin{pf}
Let $S$ be a degenerate Sklyanin algebra. Since $A$ is a Zhang twist of $A'$ the Hilbert series of 
$S$ is the same as that of 
$$
\frac{k\langle u,v,w\rangle}{(u^2,v^2,w^2)} \cong \frac{k[u]}{(u^2)} * \frac{k[v]}{(v^2)} * \frac{k[w]}{(w^2)} 
$$
where the right-hand side of this isomorphism is the free product of three copies of $k[\ve]$, the ring of dual numbers. Therefore
$$
H_S(t)^{-1}=3H_{k[\ve]}(t)^{-1} -2.
$$
Since $H_{k[\ve]}(t)=1+t$ it follows that $H_S(t)$ is as claimed.
\end{pf}

Here is another simple way to compute $H_S(t)$. 
A basis for $k\langle u,v,w \rangle/(u^2,v^2,w^2)$ is given by all words in $u$, $v$, and $w$, that do not have $uu$, $vv$, or $ww$, as a subword. The number of such words of length $n \ge 1$ is easily seen to be 
$3.2^{n-1}$. Since $3.2^{n-1}=2^n+2^{n-1}$ the Hilbert series of this algebra is $(1+t)(1-2t)^{-1}$. The same argument applies to the algebra with relations $uv=vw=wu=0$.

\section{Two quivers}
\label{sect.2.Q}

\subsection{}

The main result in \cite{HS} is the following.

\begin{thm}
\label{thm.HS}
If $A$ is a finitely presented monomial algebra there is a finite quiver $Q$ and a homomorphism $f:A \to kQ$
such that $-\otimes_A kQ$ induces an equivalence
$$
\QGr A \equiv \QGr kQ.
$$
One can take $Q$ to be the Ufnarovskii graph of $A$.
\end{thm}

A $k$-algebra is said to be {\sf matricial} if it is a finite product of matrix algebras over $k$. A 
 $k$-algebra is   {\sf ultramatricial} if it is a direct limit, equivalently a union, of matricial $k$-algebras.

\begin{thm}
\label{thm.affine}
Let $(a,b,c) \in {\frak D}$ and let $S_{a,b,c}$ be the associated 
degenerate Sklyanin algebra. There is a quiver $Q$ and an ultramatricial algebra, $R$,
both independent of $(a,b,c)$, for which there is an equivalence of categories
$$
\QGr S_{a,b,c} \equiv \QGr kQ \equiv \Mod R.
$$
\end{thm}
\begin{pf}
Since $S_{a,b,c}$ is a monomial algebra Theorem \ref{thm.HS} implies that $\QGr S_{a,b,c} \equiv \QGr kQ$
for some quiver $Q$. 
By  \cite[Thm. 1.2]{Sm1}, for every quiver $Q$ there is an ultramatricial $k$-algebra $R(Q)$
such that $\QGr kQ \equiv \Mod R(Q)$. Because $\QGr S_{a,b,c}$ is the same for all $(a,b,c) \in {\frak D}$, $Q$ and $R(Q)$ can be taken the same for all such $(a,b,c)$.
\end{pf}

\subsection{}
The Ufnarovskii graph for $k\langle u,v,w\rangle/(u^2,v^2,w^2)$ is  the quiver
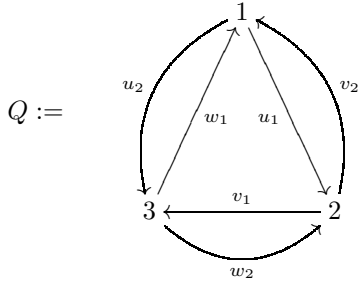
\begin{equation}
\label{Q.x2.y2.z2}
 \UseComputerModernTips
\xymatrix{ 
&& 1 \ar[ddr]_{u_1}   \ar@/_1.5pc/[ddl]_{u_2}
\\
Q:=
\\
&3 \ar[uur]_{w_1}   \ar@/_1.5pc/[rr]_{w_2} && 2 \ar[ll]_{v_1}  \ar@/_1.5pc/[uul]_{v_2}
}
\end{equation}
The equivalence $\QGr A \equiv \QGr kQ$ is induced by the $k$-algebra homomorphism 
\newline
$f:A \to kQ$ 
defined by
\begin{align*}
f(u) & = u_1+u_2, \\
f(v) & = v_1+v_2,  \\
f(w) & = w_1+w_2. 
\end{align*}

\subsection{}
The Ufnarovskii graph for $k\langle u,v,w\rangle/(uv,vw,wu)$ is  the quiver
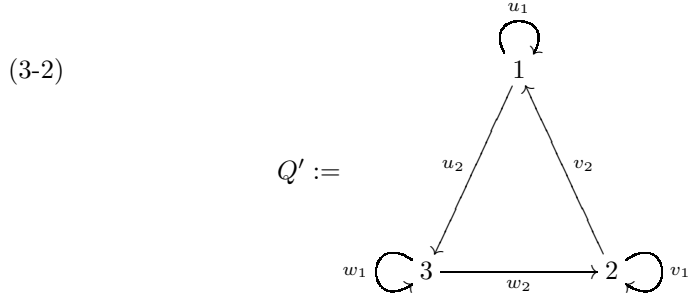
\begin{equation}
\label{Q.yx.zy.xz}
 \UseComputerModernTips
\xymatrix{ 
&&    1    \ar[ddl]_{u_2} \ar@(ul,ur)[]^{u_1}
\\
Q':= 
\\
&\ar@(ul,dl)[]_{w_1}  3   \ar[rr]_{w_2} && 2  \ar[uul]_{v_2} \ar@(ur,dr)[]^{v_1}
}
\end{equation}
The equivalence $\QGr A' \equiv \QGr kQ'$ is induced by the $k$-algebra homomorphism $f':A'\to kQ'$ 
defined by
\begin{align*}
f'(u) & = u_1+u_2, \\
f'(v) & = v_1+ v_2, \\
f'(w) & = w_1+ w_2.
\end{align*}

\subsection{Direct proof that $\QGr kQ$ is equivalent to $\QGr kQ'$}

Since $\QGr A \equiv \QGr A'$,  $\QGr kQ$ is equivalent to $\QGr kQ'$. This equivalence is not obvious 
but follows directly from \cite[Thm. 1.8]{Sm1} as we will now explain.

Given a quiver $Q$ its {\sf $n$-Veronese} is the quiver $Q^{(n)}$ that has the same vertices as $Q$ but the arrows in $Q^{(n)}$ are the paths in $Q$ of length $n$.  By \cite[Thm. 1.8]{Sm1},
$\QGr kQ \equiv \QGr kQ^{(n)}$, this being a consequence of the fact that 
$kQ^{(n)}$ is isomorphic to the $n$-Veronese subalgebra $(kQ)^{(n)}$ of $kQ$.

Hence, if the 3-Veronese quivers of $Q$ and $Q'$ are the same, then $\QGr kQ$ is equivalent to $\QGr kQ'$. 

The incidence matrix of the $n^{\th}$ Veronese quiver is the $n^{\th}$ power of the incidence matrix of the 
original quiver. The incidence matrices of $Q$ and $Q'$ are 
$$
\begin{pmatrix}
  0 & 1    & 1   \\
   1   &  0 & 1 \\
   1 & 1 & 0
\end{pmatrix}
\qquad \hbox{and} \qquad
\begin{pmatrix}
  1 & 1    & 0   \\
   0   &  1 & 1 \\
   1 & 0 & 1
\end{pmatrix}
$$
and the third power of each  is 
$$
\begin{pmatrix}
  2 & 3    & 3   \\
   3   &  2 & 3 \\
   3 & 3 & 2
\end{pmatrix}
$$
so $\QGr kQ \equiv \QGr kQ'$.

\subsection{The ultramatricial algebra $R$ in Theorem \ref{thm.affine}}

The ring $R$ in Theorem \ref{thm.affine} is only determined up to Morita equivalence but we can take it to be the 
algebra that is associated to $Q$ in \cite[Thm. 1.2]{Sm1}. That algebra has a
stationary Bratteli diagram that begins 
\begin{equation}
\label{B.diag}
 \UseComputerModernTips
\xymatrix{ 
1 \ar@{-}[drr] \ar@{-}[dr] & 1\ar@{-}[dl] \ar@{-}[dr] & 1 \ar@{-}[dl] \ar@{-}[dll]
\\
2 \ar@{-}[drr] \ar@{-}[dr] & 2\ar@{-}[dl] \ar@{-}[dr] & 2 \ar@{-}[dl] \ar@{-}[dll]
\\
4 & 4 & 4
}
\end{equation}
and so on.  More explicitly, $R = \liminj R_n$ where
$$
R_n = M_{2^n}(k) \oplus M_{2^n}(k) \oplus M_{2^n}(k)
$$
and the map $R_n \to R_{n+1}$ is given by 
$$
(r,s,t) \mapsto \Bigg( \begin{pmatrix}   s    &   0 \\    0   &   t \end{pmatrix},
\begin{pmatrix}   t   &   0 \\    0   &   r \end{pmatrix},
\begin{pmatrix}   r    &   0 \\    0   &   s \end{pmatrix}
\Bigg).
$$

\begin{prop}
The ring $R$ in Theorem \ref{thm.affine} is left and right coherent, non-noetherian, simple, 
and von Neumann regular.
\end{prop}
\begin{pf}
It is well-known that ultramatricial algebras are von Neumann regular and left and right coherent.

The simplicity of $R$ can be read off from the shape of its Bratteli diagram. 
Let $x$ be a non-zero element of $R$.
There is an integer $n$ such that $x \in R_n$ so $x=(x_1,x_2,x_3)$ where each 
$x_i$ belongs to one of the matrix factors of $R_n$. 
Some $x_i$ is non-zero so, as can be seen from the Bratteli diagram, the image of $x$ in $R_{n+2}$
has a non-zero component in each matrix factor of $R_{n+2}$. The ideal of $R_{n+2}$ generated by the image of $x$ is therefore $R_{n+2}$. Hence $RxR=R$.

The ring $R$ is not noetherian because its identity element can  be written as a sum of arbitrarily many mutually orthogonal idempotents.
\end{pf}

Because $R$ is ultramatricial it is unit regular \cite[Ch. 4]{vNR} hence directly finite \cite[Ch. 5]{vNR},
and it satisfies the comparability axiom \cite[Ch. 8]{vNR}.

\subsection{$Q$ and $Q'$ are McKay quivers}
\label{sect.McKay}

Section 2 of \cite{Sm4} describes how to associate a McKay quiver to the action of a finite abelian group acting 
semisimply on a $k$-algebra. Both $Q$ and $Q'$ can be obtained in this way.

\begin{lemma}
\label{lem.equivar.F}
Let $F=k \langle X,Y \rangle$ be the free algebra with $\deg X=\deg Y=1$.
Let $\mu_3$ be the group of $3^{\rm rd}$ roots of unity
in $k$.
\begin{enumerate}
  \item 
If $\a:\mu_3 \to \Aut_{\sf gr.alg} F$ is the homomorphism such that $\xi\hdot X = \xi X$ and  $\xi\hdot Y= \xi^2 Y$, then $kQ \cong F \rtimes_{\a} \mu_3$. 
  \item 
If $\b:\mu_3 \to \Aut_{\sf gr.alg} F$ is the homomorphism such that $\xi\hdot X = \xi X$ and  $\xi\hdot Y= Y$, then $kQ' \cong F \rtimes_{\b} \mu_3$. 
\end{enumerate}
\end{lemma}
\begin{pf}
This is a special case of \cite[Prop. 2.1]{Sm4}.
\end{pf}

Even without Lemma \ref{lem.equivar.F} one can see that $kQ$ is Morita equivalent to $F \rtimes_\a \mu_3$ and 
$kQ'$ is Morita equivalent to $F \rtimes_\b \mu_3$ because, as we will explain in the next paragraph, 
the category of $\mu_3$-equivariant $F$-modules is equivalent to the category of representations of 
$Q$ (or $Q'$, depending on the action of $\mu_3$ on $F$).

Given a $\mu_3$-equivariant $F$-module $M$ let 
$$
M_i=\{m \in M \; | \; \xi \hdot m = \xi^i m \; \hbox{for all $\xi \in \mu_3$}\}
$$
and place $M_i$ at the vertex labelled $i$ in $Q$ or $Q'$. In Lemma \ref{lem.equivar.F}(1), 
the clockwise arrows give the action of $X$ on each $M_i$ and the  counter-clockwise arrows give 
the action of $Y$ on each $M_i$. (In other words, there is a homomorphism $F \to kQ$ given by $X \mapsto
u_1+v_1+w_1$ and $Y \mapsto u_2+v_2+w_2$.) This functor from $\mu_3$-equivariant $F$-modules
to $\Mod kQ$ is an equivalence of categories.

The foregoing is ``well known'' to those to whom it is common knowledge.

\begin{prop}
Let $k$ be a field such that $\mu_3$, its $3^{\rm rd}$ roots of unity, has order 3. 
Let $F$ be the free $k$-algebra on two generators placed in degree one. Let $\a:\mu_3 \to \Aut_{\sf gr.alg} F$ be a homomorphism such that $F_1$ is a sum of two non-isomorphic representations of $\mu_3$.
If $S$ is a degenerate Sklyanin algebra, then
$$
\QGr S \equiv \QGr (F \rtimes_\a \mu_3).
$$
\end{prop}
 
 \subsubsection{Remark}
In \cite{Sm0}, it is shown that $\QGr F$ is equivalent to $\Mod T$ where $T$ is the ultramatricial algebra 
with Bratteli diagram  
$$
 \UseComputerModernTips
\xymatrix{ 
1 \ar@/^.7ex/[r]   \ar@/_.7ex/[r] & 2 \ar@/^.7ex/[r]   \ar@/_.7ex/[r] & 4 \ar@/^.7ex/[r]   \ar@/_.7ex/[r] & 8  \cdots
}
$$
We do not know an argument that explains the relation between this diagram and that in (\ref{B.diag})
which arises in connection with $\QGr (F \rtimes \mu_3)$.
 
 \subsection{The Grothendieck group of $\qgr S$}

 \begin{prop}
 Let $S$ be a degenerate Sklyanin algebra. Then
 $$
 K_0(\qgr S) \cong   \ZZ\bigg[\frac{1}{8}\bigg] \oplus \ZZ \oplus \ZZ
$$
with $[\cO]=(1,0,0)$.
\end{prop}
\begin{pf}
Since $\qgr S \equiv \mod R$,  
$$
K_0(\qgr S) \cong  K_0(R) \cong \liminj K_0(R_n)
$$  
where the second isomorphism holds becaus $K_0(-)$ commutes with direct limits.

Since $R_n$ is a product of three matrix algebras, $K_0(R_n) \cong \ZZ^3$. We can pick these isomorphisms 
in such a way that the map $K_0(R_n) \to K_0(R_{n+1})$ induced  by the inclusion $\phi_n:R_n \to R_{n+1}$ in the Bratteli diagram  is left multiplication by 
\begin{equation}
\label{M.matrix}
M:= \begin{pmatrix}
  0 &1&1   \\
    1  & 0 & 1  \\
   1  & 1 & 0 
\end{pmatrix}.
\end{equation}
The directed system  $\cdots \to K_0(R_n) \to K_0(R_{n+1}) \to \cdots$ is therefore isomorphic to 
$$
\cdots \to \ZZ^{3} \stackrel{M}{\longrightarrow}  \ZZ^{3} \stackrel{M}{\longrightarrow} \ZZ^{3} \stackrel{M}{\longrightarrow} \cdots.
$$

As noted above,
$$
M^3:= \begin{pmatrix}
  2 &3&3   \\
  3  & 2 & 3  \\
   3  & 3 & 2 
\end{pmatrix}.
$$
We have
$$
M^3   \begin{pmatrix} 1 \\ 1 \\ 1  \end{pmatrix}
= 8  \begin{pmatrix} 1 \\ 1 \\ 1 \end{pmatrix},
\qquad 
M^3   \begin{pmatrix} -1 \\ 1 \\ 0  \end{pmatrix}
= -  \begin{pmatrix} -1 \\ 1 \\ 0  \end{pmatrix},
\qquad 
M^3    \begin{pmatrix} 0 \\ -1 \\ 1   \end{pmatrix}
= -  \begin{pmatrix} 0 \\ -1 \\ 1  \end{pmatrix}.
$$
Since $\liminj K_0(R_n)$ is also the direct limit of the directed system
$$
\cdots \to \ZZ^{3} \stackrel{M^3}{\longrightarrow}  \ZZ^{3} \stackrel{M^3}{\longrightarrow} \ZZ^{3} \stackrel{M^3}{\longrightarrow} \cdots, 
$$
 $K_0(R)$ is isomorphic to $ \ZZ\big[\frac{1}{8}\big] \oplus \ZZ \oplus \ZZ$.
Together with the observation that $\cO \cong \cO(-1) \oplus \cO(-1)$ this completes the proof.
\end{pf}

\section{Point modules for $S_{a,b,c}$}

\subsection{}
A {\sf point module} over a connected graded $k$-algebra $A$ is a graded $A$-module $M=M_0\oplus M_1 \oplus
\cdots$ such that $\dim_k M_n=1$ for all $n \ge 0$ and $M=M_0A$.

\subsection{}
If one connected graded algebra is a Zhang twist of another their categories of graded modules are equivalent via a functor that sends point modules to point modules. Thus, to determine the point modules over 
a degenerate Sklyanin algebra it suffices to determine the point modules over $A$, the algebra with 
relations $u^2=v^2=w^2=0$.

\subsection{}
  
The letters $u,v,w$ will now serve double duty: they are elements in $A$ and are also an ordered set of  homogeneous coordinates on $\PP^2=\PP(A_1^*)$, the lines in $A_1^*$.

If $M=ke_0 \oplus ke_1 \oplus \cdots$ is a point module with $\deg e_n=n$, we define points $p_n \in \PP^2$,
$n \ge 0$, by 
$$
p_n=(\a_n,\b_n,\gamma_n)
$$
where $e_n\hdot u=\a_n e_{n+1}$,  $e_n\hdot v=\b_n e_{n+1}$,  and $e_n\hdot w=\c_n e_{n+1}$. 
The $p_n$s do not depend on the choice of homogeneous basis for $M$. We call $p_0,p_1,\ldots$ the {\sf point sequence} associated to $M$.

\subsection{}
Let $E$ be the three lines in $\PP^2$ where $uvw=0$.  
We call the points that lie on two of those lines 
{\it intersection points}. If $p$ is an intersection point the component of $E$ that does not pass through $p$
is called {\it the line opposite $p$}. If $L$ is a component of $E$ we call the intersection 
point that does not lie on $L$ {\it the point opposite $L$}.

\begin{lemma}
\label{lem.pt.mod.seq}
Let $(p_0,p_1,\ldots)$ be the point sequence associated to a point module $M$.
\begin{enumerate}
\item{}
Every $p_n$ lies on $E$.
  \item 
  If $p_n$ is an intersection point
then $p_{n+1}$ lies on the line opposite $p_n$ and can be any point on that line. 
  \item 
  If $p_n$ is not an intersection point and lies on $L$, 
then $p_{n+1}$ is the point opposite $L$.   
\end{enumerate} 
\end{lemma}
\begin{pf}
(1)
If $p_n \notin E$, then $e_n\hdot u$, $e_n\hdot v$, and $e_n\hdot w$ are non-zero scalar multiples of 
$e_{n+1}$. But $u^2$, $v^2$, and $w^2$, are zero in $A$ so $e_n\hdot u^2=e_n\hdot v^2=e_n\hdot w^2=0$
which implies that $e_{n+1}A=ke_{n+1}$ in contradiction of the fact that $M$ is a point module.

(2) 
Suppose $p_n$ lies on the lines $u=0$ and $v=0$.  The line opposite $p_n$ is the line $w=0$.
Since $e_n\hdot u=e_n\hdot v=0$, $e_n\hdot w$ must be a non-zero multiple of $e_{n+1}$ so,
since $w^2=0$, $e_{n+1}\hdot w=0$; i.e., $w(p_{n+1}) \ne 0$. The other cases are similar.

(3)
Suppose $p_n$ lies on the line $u=0$ but not on $v=0$ or $w=0$. Then $e_n \hdot v=e_n\hdot w = e_{n+1}$
and, because $v^2=w^2=0$,  $e_{n+1} \hdot v=e_{n+1}\hdot w =0$. Hence $p_{n+1}$ is the point opposite the line $u=0$. The other cases are similar.
\end{pf}

Lemma \ref{lem.pt.mod.seq} shows that \cite[Thm. 1.7]{W} is not correct. 

\subsection{}
 If $s$ is a non-zero word of length $n$ we write $s^\perp$ for the other non-zero words of length $n$.

A sequence $p_0,p_1,\ldots,p_n$ of points in $E$ that arise from a (truncated) point module of length $n+2$ 
is said to be {\sf special} if each $p_i$ is an intersection point. A (truncated) point module giving rise to a special
sequence is also called {\sf special}. 

\begin{lem}
If $n \ge 0$, there are $2^n\times 3$ special point sequences $p_0,\ldots,p_n$ and $2^n\times 3$ special truncated point modules of length $n+2$ up to isomorphism. 
\end{lem}

As noted before, $\dim A_{n+1} =2^n \times 3$ for all $n \ge 1$. 

\begin{prop}
\label{prop.sp.pt.mods}
If $s=s_0\ldots s_n$ is a non-zero word of length $n+1$, there is a special truncated point module $N$ of length $n+2$ such that 
$N \hdot s \ne 0$ and $N \hdot t=0$ for every $t \in s^\perp$.
\end{prop}
\begin{pf}
We make $N:=ke_0 \oplus \cdots \oplus ke_{n+1}$ into a right $A$-module by having $a \in \{u,v,w\}$ act 
as follows:
$$
e_i \hdot a = 
\begin{cases} 
	e_{i+1} & \text{ if $a=s_{i}$}
	\\
	0 & \text{ if $a \ne s_{i}$}
\end{cases}
$$
for $0 \le i \le n$ and $e_n\hdot u = e_n\hdot v = e_n\hdot w =0$. Because $s_i \ne s_{i+1}$, $u^2$, $v^2$, and $w^2$, act as zero on $N$. By construction, $e_0 \hdot s = e_{n+1} \ne 0$ so $N$ is a truncated point module. If $t=t_0\ldots t_n \in s^\perp$, then $t_i \ne s_i$ for some $i$, whence $e_0\hdot t=0$. Of course, $e_i \hdot t=0$ for $i \ge 1$ so $N\hdot t=0$.   
 
The point sequence associated to $N$, $p_0,\ldots,p_n$, is given by
$$
p_i = 
\begin{cases}
(1,0,0) & \text{if $s_i=u$}
\\
(0,1,0) & \text{if $s_i=v$}
\\
(0,0,1) & \text{if $s_i=w$.}
\end{cases}
$$
In particular, this is a special sequence and $N$ is therefore a special truncated point module.
\end{pf} 

\begin{prop}
\label{prop.trunc.pt.mod}
Let $a \in A_n-\{0\}$.
There is a truncated special point module $N$ of length $n+1$ such that $N \hdot a \ne 0$.
\end{prop}
\begin{pf}
Let 
$$
I = \{b \in A \; | \; Nb=0 \; \hbox{ for all truncated point modules $N$ of length $n+1$}\}.
$$
To prove the proposition we must show that $I_n =0$.

Let $s$ be a non-zero word of length $n$. By Proposition \ref{prop.sp.pt.mods}, there is a truncated point module $N$ of length $n+1$ such that $(\Ann N)_n \subset s^\perp$. Hence $I_n \subset \Span(s^\perp)$. However, $L_n$ is a basis for $A_n$ so the intersection of $\Span(s^\perp)$ as $s$ ranges over $L_n$
is zero. Hence $I_n=0$.
\end{pf}

\begin{cor}
If $a$ is a non-zero element in $A$ there is a special point module $M$ such that $M\hdot a \ne 0$.
\end{cor}

Proposition \ref{prop.trunc.pt.mod} implies that the natural map from $A$ to its point parameter ring (see 
\cite[Sect. 1]{W}) is injective.
Hence \cite[Thm. 1.9]{W} and \cite[Cor. 1.10]{W} are incorrect.

 \end{document}